\documentclass[12pt]{article}
\usepackage{amsmath}
\usepackage{amssymb}
\usepackage{vatola}

\def\q{\quad}

\def\t{\text}
\def\qtq#1{\q\t{#1}\q}
\def\f{\frac}

\def\a{\alpha}

\def\kn1{k_1!\cdots k_n!}

\def\b{\binom}

\def\ep{\varepsilon}

\def\({\left(}
\def\){\right)}

\let \pro=\proclaim
\let \endpro=\endproclaim

\begin{document}

 \centerline {\bf Tur\'an's problem and generalized Ramsey numbers}
$$\q$$
\centerline{Zhi-Hong Sun} $$\q$$ \centerline{School of Mathematical
Sciences, Huaiyin Normal University,} \centerline{Huaian, Jiangsu
223001, P.R. China} \centerline{E-mail: zhihongsun@yahoo.com}
\centerline{Homepage: http://www.hytc.edu.cn/xsjl/szh}

\abstract{Let $n,r,k,s$ be positive integers with $n,k\ge 2$. The
generalized Ramsey number $R(n,r;k,s)$ is the smallest positive
integer $p$ such that for every graph $G$ of order $p$, either $G$
contains a  subgraph induced by $n$ vertices with at most $r-1$
edges, or the complement $\overline G$ of $G$ contains a subgraph
induced by $k$ vertices with at most $s-1$ edges. In this paper  we
completely determine $R(n,n(n-1)/2-r;k,1)$ for $n\ge 4$ and $r\le
n-2$, and pose several conjectures on Ramsey numbers.
\par\q
}
 \endabstract
 \footnotetext[1] {The author is
supported by the National Natural Science Foundation of China (grant
No. 11371163).}

\section*{1. Introduction}
\par Let $n,k\ge 2$ be positive
integers. The classical Ramsey number $R(n,k)$ is the minimum
positive integer such that every graph on $R(n,k)$ vertices has a
complete subgraph $K_n$ or an independent set with $k$ vertices. Up
to now we only know the following exact values of Ramsey numbers
(see [1, p.187] and [8, p.4]):
$$\aligned &R(3,3)=6,\  R(3,4)=9,\  R(3,5)=14,\  R(3,6)=18,\
R(3,7)=23,\\& R(3,8)=28,\ R(3,9)=36,\ R(4,4)=18,\ R(4,5)=25.
\endaligned\tag 1.1$$
\par
 In this paper
we introduce the following  generalized Ramsey numbers.
\pro{Definition 1.1} Let $n,r,k,s$ be positive integers with $n,k\ge
2$. We define the generalized Ramsey number $R(n,r;k,s)$ to be the
smallest positive integer $p$ such that for every graph $G$ of order
$p$, either $G$ contains a subgraph induced by $n$ vertices with at
most $r-1$ edges, or the complement $\overline G$ of $G$ contains a
 subgraph induced by $k$ vertices with at most $s-1$ edges.\endpro
\par Clearly
$R(n,1;k,1)=R(n,k)$. In 1981, Bolze and Harborth [2] introduced the
generalized Ramsey number $r_{m,n}(s,t)=R(m,\b m2-s+1;n,\b n2-t+1)$
$(1\le s\le \b m2,1\le t\le \b n2)$. For some values of
$r_{m,n}(s,t)$, see [2,5,7].

\par Let $L$ be a family of forbidden graphs. As usual we use $ex(p;L)$
 to denote the maximal number of edges in a graph of order $p$ excluding
  any graphs in $L$. The general Tur\'{a}n 's problem is to
 evaluate $ex(p;L)$.
 \par
 For two positive integers $p$ and
 $k$, let $p_0$ be the least nonnegative residue of $p$ modulo $k$, and
 let
 $$ t_k(p)=\Big(1-\f 1k\Big)\f{p^2-p_0^2}2+\f{p_0^2-p_0}2.\tag 1.2$$
 In 1941 Tur\'an [11] showed that
 $ex(p;K_k)=t_{k-1}(p)$ for $k>1$, where  $K_k$ is the complete graph with $k$
vertices. This is now called Tur\'an's theorem. In 1963, Dirac [4]
proved a vast generalization of Tur\'an's theorem, see (2.6). Let
$p,n$ and $t$ be positive integers satisfying $2\leq t\leq\f n2+2$
and $p\geq n\geq 4$. With the help of Dirac's generalization of
Tur\'an's theorem, in Section 3 we show that if $G$ is a graph of
order $p$ and every induced subgraph of $G$ by $n$ vertices has at
most $n-t$ edges, then
$$\a(G)\geq \Big[\f{p-[\f{n+4-2t}3]}2\Big]+1,\tag 1.3$$
where  $\alpha(G)$ is the independence number of $G$ and $[x]$ is
the greatest integer not exceeding $x$. In Section 4 we show that
$$R(n,n(n-1)/2-r;k,1)=\cases max\{n,k+r\}&\t{if}\q r\le \f
n2-1,\\ max\{n,2k-2+[\f{2r+4-n}3]\} &\t{if}\q r> \f n2-1,\endcases
\tag 1.4$$ where $k,n,r$ are positive integers with $k\ge 2$, $n\ge
4$ and $r\leq n-2$.  In the special case $k=n$, (1.4) is known, see
[5, Theorem 4].  In Section 5 we obtain a upper bound for
$R(4,3;k,1)$, and in Section 6 we pose several conjectures on
$R(n,k)$ and $R(n,r;k,s)$.

\par In addition to the above notation, all graphs in the paper are
simple graphs, throughout the paper we use the following notation:
\par $\Bbb N\f{\q}{\q}$the set of positive integers,
$\{x\}\f{\q}{\q}$the fractional part of $x$, $V(G)\f{\q}{\q}$the set
of vertices in the graph $G$, $e(G)\f{\q}{\q}$the number of edges in
the graph $G$, $d(v)\f{\q}{\q}$the degree of the vertex $v$ in a
graph, $\Gamma(v)\f{\q}{\q}$the set of those vertices adjacent to
the vertex $v$,  $\delta(G)\f{\q}{\q}$the minimal degree of $G$,
$\Delta(G)\f{\q}{\q}$the maximal degree of $G$, $g(G)\f{\q}{\q}$the
girth of $G$,
 $C_n\f{\q}{\q}$the
cycle with $n$ vertices,   $G[V]\f{\q}{\q}$the subgraph of $G$
induced by the vertices in $V$, $G-V\f{\q}{\q}$the subgraph of $G$
obtained  by deleting the vertices in $V$ and all edges incident
with them.

 \section*{2. Some applications of the generalization of
Tur\'an's theorem}
\par Let $L$ be a family of forbidden graphs, and let $p\in\Bbb N$.
The famous Erd\"{o}s-Stone Theorem states
 that(see [1, pp.122-123])
$$ex(p;L)=\f 12\Big(1-\f 1{\chi(L)-1}\Big)p^2+o(p^2),\tag 2.1
$$where $\chi (L)=\min \{\chi(G):G\in L\}$ and $\chi (G)$ is the
chromatic number of $G$.

\pro{Lemma 2.1} Let $L$ be a family of forbidden graphs. Let $G$
be a graph of order $p\ge 2$ without any graphs in $L$ and
$e(G)=ex(p;L)$. Then
$$\delta(G)\ge  ex(p;L)-ex(p-1;L).$$
\endpro
Proof. Suppose that $v$ is a vertex in $G$ such that
$d(v)=\delta(G)$. Then clearly
$$ex(p;L)-\delta(G)=e(G)-d(v)=e(G-v)\le ex(p-1;L),$$
which yields the result.
 \pro{Theorem 2.1} Let $L$ be a family of
forbidden graphs. Let $p>2$ be a positive integer. Then
$$ ex(p;L)\leq \Big[\f {p\cdot ex(p-1;L)}{p-2}\Big]. $$
\endpro
Proof. Suppose that $G$ is a graph of order $p$ without any graphs
in $L$ and $e(G)=ex(p;L)$. Using Lemma 2.1 and Euler's theorem we
see that
 $$p\big(ex(p;L)-ex(p-1;L)\big)\le p\cdot \delta(G)\le \sum_{v\in V(G)}d(v)
=2\;e(G)=2\;ex(p;L).$$
  Thus, $$ ex(p;L)\leq \f {p\cdot ex(p-1;L)}{p-2}. $$
  As $ex(p;L)$ is an integer, the result follows.
\pro{Definition 2.1} Let $n$ and $m$ be two integers such that
$n\ge 1$ and
 $m\ge 0$. A graph $G$ is said to be an $(n,m)$ graph
if every induced subgraph by $n$ vertices has at most $m$ edges.
For $p\ge n$ we define $e(n,m;p)$ to be the maximal number of
edges in $(n,m)$ graphs of order $p$.\endpro
\par By Theorem 2.1, we have
$$ e(n,m;p)\leq \Big[\f {p\cdot e(n,m;p-1)}{p-2}\Big].\tag 2.2 $$
This has been given by the author in [9].
\par Let $p,m,n,r\in\Bbb N$ with $p\ge m\ge n\ge 3$ and $r\le \b
n2$. Let $G$ be an $(n,r)$ graph of order $p$ with $e(G)=e(n,r;p)$.
Then clearly $G$ is  an $(m,e(n,r;m))$ graph of order $p$. As
$e(G)=e(n,r;p)$, we have
$$e(n,r;p)\le e(m,e(n,r;m);p).\tag 2.3$$
 \par In 1991, the author[9] proved that for $p\ge n\ge 3$,
$$e(n,n-2;p)=\Big[\f{(n-2)p}{n-1}\Big]
\qtq{and}e(n,n-1;p)=ex(p;\{C_3,\ldots, C_n\}).\tag 2.4$$ The first
formula in (2.4) was first proved by Gol'berg and Gurvich in [6].
 Let
$t_k(p)$ be given by (1.2), and let $k\geq 2$ be the unique integer
such that $t_{k-1}(n)\leq m<t_k(n)$. Using the Erd\"{o}s-Stone
Theorem, in 1991 the author proved that (see [9])
$$e(n,m;p)\sim \f {k-2}{2(k-1)}p^2 \ \ as \ \ p \rightarrow
+\infty.\tag 2.5$$
 In [4] Dirac
extended Tur\'{a}n's theorem by proving the
 following result (see also [9]):
 $$e(n,t_{k-1}(n);p)=t_{k-1}(p)\qtq{for}p\ge n\ge k\ge 2.\tag 2.6$$
 For
$n\ge 2m\ge 2$ it is clear that $t_{n-m}(n)=\b n2-m$. Thus,
$$e\big(n,\f{n(n-1)}2-m;p\big)=t_{n-m}(p)\qtq{for}p\ge n\ge 2m\ge 2.\tag
2.7$$

  \pro{Theorem 2.2} Let $p,n,k\in\Bbb N$ with $p\ge n>k$, and let
   G be a graph of order $p$. If
  $$e(G)<\Big[\f pk\Big]p-\f{[\f pk]([\f pk]+1)}2k,$$
  then $G$ contains a subgraph induced by $n$ vertices
  with at most $[\f nk](n-\f{[\f nk]+1}2k)-1$ edges.
 \endpro
 Proof. Set $s=[\f pk]$. Clearly $p=ks+r$ for some $r\in \{ 0,1,\cdots ,k-1\}.$
 Thus
 $$\aligned \b p2-t_k(p)&=\b p2-\f
 {k-1}k\cdot\f{p^2-r^2}2-\f{r^2-r}2\\&=\f{p-r}
2\big(\f{p+r}k-1\big)=\f{ks}2\big(\f{p+p-ks}k-1\big)
\\&=sp-\f{s(s+1)}2k.\endaligned$$
Hence, if $e(G)<sp-\f{s(s+1)}2k$, then $e(G)<\b p2-t_k(p)$ and so
$e(\overline G)=\b p2-e(G)>t_k(p)=e(n,t_k(n);p)$ by (2.6).
Therefore, $\overline G$ contains an induced subgraph by $n$
vertices with at least $t_k(n)+1$ edges. Hence $G$ contains an
induced subgraph by $n$ vertices with at most  $\b n2-t_k(n)-1$
edges. As $\b n2-t_k(n)=[\f nk](n-\f{[\f nk]+1}2k)$, we deduce the
result.
 \pro{Corollary 2.1} Let $p,n,k\in\Bbb N$ with $2\le k+1\le n\le 2k$ and $p\ge n$.
   Let G be a graph of order $p$ satisfying
   $$e(G)<\Big[\f pk\Big]p-\f{[\f pk]([\f pk]+1)}2k.$$
   Then $G$ contains a subgraph induced by $n$ vertices
  with at most $n-k-1$ edges. In particular, for $p\ge n=k+1$ we have $\alpha(G)\ge k+1$.
 \endpro
 Proof. Observe that $[\f nk]=1$ or $2$ according as
 $k+1\le n<2k$ or $n=2k$. We then have $[\f nk](n-\f{[\f
 nk]+1}2k)-1=n-k-1$. Now applying Theorem 2.2 we deduce the result.
\pro{Corollary 2.2} Let G be a graph of order $p\geq 4$. Let $t$ be
a nonnegative integer with $t< \f p4-3\{\f{p-t}3\}$. If $e(G)<
p+2t+6\{\f{p-t}3\}$, then $\alpha(G)\geq [\f{p-t}3]+1$.
\endpro
Proof. Let $r=p-3[\f{p-t}3].$ Then
$r=p-3(\f{p-t}3-\{\f{p-t}3\})=t+3\{\f{p-t}3\}\geq 0$. Since $t<\f
p4-3\{\f{p-t}3\}$ we have $4\{\f{p-t}3\}<\f{p-4t}3$ and so
$$r=t+3\big\{\f{p-t}3\big\}<\f{p-t}3-\{\f{p-t}3\}=\big[\f{p-t}3\big].$$
 Thus $[\f p{[\f{p-t}3]}]=3.$
 Hence,
 $$\align &\Big[\f p{[\f{p-t}3]}\Big]\Big(p-\f{\big[\f
p{[\f{p-t}3]}\big]+1}2\big[\f{p-t}3\big]\Big)\\&=
3\big(p-2\big[\f{p-t}3\big]\big)
=3p-6\big(\f{p-t}3-\big\{\f{p-t}3\big\}\big)
=p+2t+6\big\{\f{p-t}3\big\}.
\endalign$$
Since $t\ge 0>3-2p$ we see that $p>\f{p-t}3+1\ge [\f{p-t}3]+1$. We
also have $[\f{p-t}3]\ge [\f{p-p/4}3]=[\f p4]\ge 1$.
 Now applying the above and taking $k=[\f{p-t}3]$ in
Corollary 2.1 we obtain the result.
 \pro{Theorem 2.3}
Let $G$ be a graph of order $p\geq 4$. Let $t$ be a nonnegative
integer so that $3\{\f{p+t}3\}-\f p4<t\leq \f p2+3\{\f{p+t}3\}$. If
$$e(G)<p-2t+3\Big\{\f{p+t}3\Big\}+\t{max}\ \Big\{t,3\Big\{\f{p+t}3\Big\}
\Big\},$$
 then $\alpha(G)\geq
[\f{p+t}3]+1.$
\endpro
Proof. It is clear that
$$\align &p-2\big[\f{p+t}3\big]=p-2\big(\f{p+t}3-\big\{\f{p+t}3\big\}\big)=\f
13\big(p-2t+6\big\{\f{p+t}3\big\}\big)\geq
0,\\&p-3\big[\f{p+t}3\big]=p-3\big(\f{p+t}3-\big\{\f{p+t}3\big\}\big)
=3\big\{\f{p+t}3\big\}-t\endalign$$ and
$$p-4\big[\f{p+t}3\big]=p-4\big(\f{p+t}3-\big\{\f{p+t}3\big\}\big)
=-\f{p+4t-12\{\f{p+t}3\}}3<0.$$
 Thus,
 $$\Big[\f p{[\f{p+t}3]}\Big]=\cases 2
&\t{if}\q t>3\{\f{p+t}3\},\\3 & \t{if}\q
t\leq3\{\f{p+t}3\}.\endcases$$ Set $k=[\f{p+t}3]$. We see that
$$\aligned &\big[\f pk\big]\Big(p-\f{[\f pk]+1}2k\Big)\\&=\cases
2p-3k=2p-3(\f{p+t}3-\{\f{p+t}3\})=p-t+3\{\f{p+t}3\} &\t{if}\q
t>3\{\f{p+t}3\},\\3p-6k=3p-6(\f{p+t}3-\{\f{p+t}3\})=p-2t+6\{\f{p+t}3\}&
\t{if}\q t\leq3\{\f{p+t}3\},\endcases
\\&=p-2t+3\Big\{\f{p+t}3\Big\}+\t{max}\ \Big\{t,3\Big\{\f{p+t}3\Big\}
\Big\}.\endaligned$$ Since $p\ge 4$ and $0\le t<\f p2+3$ we see that
$[\f{p+t}3]\ge 1$ and $p\ge \f{p+\f p2+3}3+1\ge [\f{p+t}3]+1$. Now
applying the above and Corollary 2.1 we obtain the result.
\par Putting $t=0,1,2$ in Theorem 2.3 we deduce the following
corollaries.
 \pro{Corollary 2.3} Let $G$ be a graph of
order $p\geq 9$. If $e(G)<p+6\{\f p3\}$, then $\alpha (G)\geq [\f
p3]+1$.\endpro
 \pro{Corollary 2.4} Let $G$ be a
graph of order $p\geq 5$. If
$$\aligned e(G)<\cases p &\text{if } 3\mid p,\\p+2 &\text{if }3\mid p-1,\\p-1 &\text{if }
3\mid p-2,\endcases\endaligned$$ then $\alpha (G)\geq
[\f{p+1}3]+1$.\endpro
 \pro{Corollary 2.5}Let $G$ be a graph of
order $p\geq 4$. If $$\aligned e(G)<\cases p &\text{if } 3\mid
p,\\p-2 &\text{if }3\mid p-1,\\p-1 &\text{if } 3\mid
p-2,\endcases\endaligned$$ then $\alpha (G)\geq
[\f{p+2}3]+1$.\endpro
 \pro{Theorem 2.4} Let $p,n,t\in\Bbb N$ with $t\leq \f
 n2+2$ and $p>n+7-2t.$ If $G$ is a graph of order $p$ and
  $e(G)<\f {p+n}2+1-t,$ then $\a(G)\geq\big[\f{p-[\f{n+4-2t}3]}2\big]+1.$
  \endpro
  Proof. Set $r=p-2\big[\f{p-[\f{n+4-2t}3]}2\big].$ Then
  $$ r\geq p-2\cdot\f{p-[\f{n+4-2t}3]}2=\big[\f{n+4-2t}3\big]\geq 0$$
  and
 $$ r=p-2\Big(\f{p-[\f{n+4-2t}3]}2-\Big\{\f{p-[\f{n+4-2t}3]}2\Big\}\Big)
=2\Big\{\f{p-[\f{n+4-2t}3]}2\Big\}+\Big[\f{n+4-2t}3\Big].$$ Hence
$$\align &\Big[\f{p-[\f{n+4-2t}3]}2\Big]-r
\\&=\f{p-[\f{n+4-2t}3]}2-\Big\{\f{p-[\f{n+4-2t}3]}2\Big\}
-2\Big\{\f{p-[\f{n+4-2t}3]}2\Big\}-[\f{n+4-2t}3]
\\&=\f{p-3[\f{n+4-2t}3]}2-3\Big\{\f{p-[\f{n+4-2t}3]}2\Big\}
\\&=\f{p-3\big(\f{n+4-2t}3-\big\{\f{n+4-2t}3\big\}\big)}2
-3\Big\{\f{p-[\f{n+4-2t}3]}2\Big\}
\\&=\f{p-(n+4-2t)}2+\f 32\Big\{\f{n+4-2t}3\Big\}
-3\Big\{\f{p-[\f{n+4-2t}3]}2\Big\} >\f 32+0-3\cdot\f 12=0.
\endalign$$
So we have
$$\bigg[\f p{[\f{p-[\f{n+4-2t}3]}2]}\bigg]=2.$$
As $$\align  2p-\f{2\times 3}2\Big[\f{p-[\f{n+4-2t}3]}2\Big]
&=2p-3\Big(\f{p-[\f{n+4-2t}3]}2-\Big\{\f
{p-[\f{n+4-2t}3]}2\Big\}\Big)
\\&=\f p2+\f 32\Big(\f{n+4-2t}3-\Big\{\f{n+4-2t}3\Big\}\Big)
+3\Big\{\f{p-[\f{n+4-2t}3]}2\Big\}
\\&=\f p2+\f{n+4-2t}2+3\Big\{\f{p-[\f{n+4-2t}3]}2\Big\}-\f 32
\Big\{\f{n+4-2t}3\Big\}
\\&\geq \f{p+n}2+2-t+0-1,
\endalign$$
we see that
$$ e(G)<\f{p+n}2+1-t\leq 2p-\f{2\times
3}2\Big[\f{p-[\f{n+4-2t}3]}2\Big].$$ Since $[\f{p-[\f{n+4-2t}3]}2]
\ge [\f{p-\f{n+4-2t}3}2]=[\f{3p-(n+4-2t)}6]\ge 1$ and $p\ge 2
[\f{p-[\f{n+4-2t}3]}2]$, by the above and Corollary 2.1 we deduce
the result.
 \section*{3. The best lower bound for
  independence numbers of $(n,m)$ graphs $(m\leq n-2)$}

\pro{Theorem 3.1} Let $p,m,n\in \Bbb N$ with $p\geq n\geq 4$ and
$m\le \f n2-1$. If $G$ is an $(n,m)$ graph of order $p$, then
$\alpha(G)\geq p-m.$\endpro Proof. We prove the theorem by induction
on $p$. Set $k=n-1-m$. Then $3\le k+1\le n\le 2k$ and so $[\f nk]=1$
or $2$. Thus, if $G$ is a graph of order $n$ with $e(G)\leq m$, then
$$e(G)\leq m<\Big[\f nk\Big]n-\f{[\f nk]([\f nk]+1)}2k.$$
Since $n\ge k+1$, applying Corollary 2.1 we see that $\alpha(G) \ge
k+1=n-m.$ So the result holds for $p=n$.
\par Now
assume $p>n$. Let $G$ be an $(n,m)$ graph of order $p$. We assert
that $G$ has an isolated vertex. Otherwise, we have $$e(G)=\f
12\sum_{v\in V(G)}d(v)\geq \f p2>\f n2\geq \Big[\f n2\Big].$$ We
choose $[\f n2]$ edges in $G$ and consider a subgraph $H$ induced by
the $[\f n2]$ edges. Clearly $H$ has at most $n$ vertices and
$e(H)=[\f n2]>\f n2-1\geq m.$ This contradicts the fact that $G$ is
an $(n,m)$ graph. So $G$ has an isolated vertex. Assume that $v$ is
an isolated vertex of $G$ and the result holds for all $(n,m)$
graphs of order $p-1$. By the above we have
$$\alpha(G)=1+\alpha(G-v)\geq 1+p-1-m=p-m.$$
Hence the theorem is proved by induction.

 \pro{Lemma 3.1} Let
$n,m\in\Bbb N$ with $m\leq n-2$. If $G$ is an $(n,m)$ graph of order
$n+1$, then $e(G)\leq m+1.$\endpro Proof. By (2.2) we have
$$e(n,m;n+1)\leq
\big[\f{(n+1)e(n,m;n)}{n-1}\big]=\big[\f{(n+1)m}{n-1}\big]=m+\big[\f{2m}{n-1}
\big]$$ As $m<n-1$ we see that $[\f{2m}{n-1}]<2$. Thus $$e(G)\leq
e(n,m;n+1)\leq m+\big[\f{2m}{n-1}\big]\leq m+1.$$ This is the
result.
 \pro{Lemma 3.2} Let $p,m,n\in\Bbb N$ with $p\geq n\geq
m+2.$ Let $G$ be an $(n,m)$ graph of order $p$. Then one of the
components of $G$ is a tree.\endpro
 Proof. Suppose that
$G_1,G_2,\ldots,G_r$ are all components of $G$ and $|V(G_i)|$ $=p_i$
for $i=1,2,\ldots,r$. Then $p_1+\cdots+p_r=p\geq n.$ Since $m<n-1$
and G is an $(n,m)$ graph, $G$ cannot contain a tree on $n$ vertices
as a subgraph. As every connected graph has a spanning tree, we must
have $p_i<n$ for all $i\in\{1,2,\ldots,r\}.$ Let us choose
$j\in\{1,2,\ldots,r-1\}$ so that $p_1+\cdots+p_j<n$ and
$p_1+\cdots+p_j+p_{j+1}\geq n.$ Then clearly $G_{j+1}$ has a subtree
$T$ on $n-(p_1+\cdots+p_j)$ vertices. All the $n$ vertices in
$G_1\cup\cdots\cup G_j\cup T$ induce a subgraph of $G$ with at least
$e(G_1)+\cdots+e(G_j)+e(T)$ edges. If $G_1,\cdots,G_r$ are not
trees, then $e(G_i)\geq p_i$ for $i=1,2,\ldots,r$ and so
$$e(G_1)+\cdots+e(G_j)+e(T)\ge p_1+\cdots+p_j+n-(p_1+\cdots+p_j)-1=n-1.$$
 Thus, the $n$
vertices in $G_1\cup\cdots\cup G_j\cup T$ induce a subgraph of $G$
with at least $n-1$ edges. This contradicts the assumption. Hence
one of $G_1,\cdots,G_r$ is a tree. This proves the lemma. \pro{Lemma
3.3} Let $T$ be a tree with at least $3$ vertices. Then there are
two vertices $u$ and $v$ in $T$ such that
 $ \alpha(T)>\alpha(T-\{u,v\})$.\endpro
 Proof. It is well known that $T$ has a vertex $u$ such that
 $d(u)=1.$ Let $v$ be the unique vertex of $G$ adjacent to $u$.
 Let $S$ be an independent set of $T-\{u,v\}.$ Then clearly $S\cup
 \{u\}$ is an independent set of $T$. Thus $\alpha(T)>\alpha(T-\{u,v\}).$
\pro {Theorem 3.2} Let $p,n,t\in\Bbb{N},2\leq t\leq\f n2+2$ and
$p\geq n\geq 4.$ If $G$ is an $(n,n-t)$ graph of order $p$, then
$$\a(G)\geq \Big[\f{p-[\f{n+4-2t}3]}2\Big]+1.$$
\endpro
Proof. We prove the theorem by induction on $p$. Now suppose $p=n$
or $n+1$. Let $G$ be an $(n,n-t)$ graph of order $p$. For $t\ge 4$
and $p=n$ we have $e(G)\le n-t<\f{p+n}2+1-t$. For $t\ge 4$ and
$p=n+1$, by Lemma 3.1 we also have $e(G)\le n-t+1<\f{p+n}2+1-t$.
Thus applying Theorem 2.4 we see that $\a(G)\geq
[\f{p-[\f{n+4-2t}3]}2]+1$. For $t=3$ and $p\in\{n,n+1\}$ it is
easily seen that $\big[\f{p-[\f{n-2}3]}2\big]=[\f{p+2}3]$. By Lemma
3.2 we have $e(G)\le  n-t+1=p-3$ for $p=n+1$ and $e(G)\le n-t=p-3$
for $p=n$. Thus, by Corollary 2.5 we have
$$\alpha(G)\ge\Big[\f{p+2}3\Big]+1=\Big[\f{p-[\f{n+4-2t}3]}2\Big]+1.$$
For $t=2$ and $p\in\{n,n+1\}$ it is easily seen that $\big[\f{p-[\f
n3]}2\big]=[\f{p+1}3]$. By Lemma 3.1 we have $e(G)\le  n-t+1=p-2$
for $p=n+1$ and $e(G)\le n-t=p-2$ for $p=n$. Thus, by Corollary 2.4
we have
$$\alpha(G)\ge\Big[\f{p+1}3\Big]+1=\Big[\f{p-[\f{n+4-2t}3]}2\Big]+1\qtq{for}
p\ge 5.$$ When $p=n=4$ and $t=2$, we also have $\alpha(G)\ge
2=[\f{p-[\f{n+4-2t}3]}2]+1.$ Summarizing the above we see that the
result is true for $p\in\{n,n+1\}$.
\par Now we assume $p\ge n+2$ and the result holds for all $(n,n-t)$ graphs
with at most $p-1$ vertices. Let $G$ be an $(n,n-t)$ graph of order
$p$. By Lemma 3.2, one of the components of $G$ is a tree. Let $T$
be such a tree. If $T\cong K_1$, for any vertex $u$ in $G-T$ we have
$\alpha (G)=1+\alpha(G-T)\ge 1+\alpha(G-\{T,u\})$. If $T\cong K_2$,
 it is clear that $\alpha (G)>\alpha(G-T)$. If
$T$ is a tree with at least three vertices, by Lemma 3.3 there are
two vertices $u$ and $v$ in $T$ such that $\alpha
(G)>\alpha(G-\{u,v\})$. Hence, in any cases, there are two vertices
$u$ and $v$ in $G$ so that $\alpha (G)>\alpha(G-\{u,v\})$. Clearly
$G-\{u,v\}$ is an $(n,n-t)$ graph of order $p-2$. Thus, by the
inductive hypothesis we have
$$\alpha(G)\ge 1+\a(G-\{u,v\})\ge 1+\Big[\f{p-2-[\f{n+4-2t}3]}2\Big]+1
=\Big[\f{p-[\f{n+4-2t}3]}2\Big]+1.$$
 So the theorem is proved by
induction. \section*{4. Evaluation of $R(n,n(n-1)/2-r;k,1)\ (r\le
n-2)$}
 \par Let $n,r,k$ be positive integers with $r\le\b n2$. By Definition 1.1,
  $R(n,n(n-1)/2-r;k,1)$ is the
smallest positive integer $p$ such that for any  $(n,r)$ graph $G$
of order $p$, we have $\alpha(G)\ge k$.
\par Now we are in a position to prove the following main result.
 \pro
{Theorem 4.1} Let $k,n,r\in\Bbb{N}$ with $k\ge 2$, $n\ge 4$ and
$r\leq n-2$.
 Then
$$R(n,n(n-1)/2-r;k,1)=\cases max\{n,k+r\}&\t{if}\q r\le \f
n2-1,\\ max\{n,2k-2+[\f{2r+4-n}3]\} &\t{if}\q r>\f n2-1.\endcases
$$
\endpro
Proof. We first assume $r\le \f n2-1$.  If $k+r\le n$, putting $p=n$
and $m=r$ in Theorem 3.1 we see that $\alpha(G)\ge n-r\ge k$ for any
graph $G$ of order $n$ with $e(G)\le r$. Thus, $R(n,\b n2-r;k,1)\le
n$ and so $R(n,\b n2-r;k,1)=n=\t{max}\{n,k+r\}$. Now suppose
$k+r>n$. Then $k>n-r\ge n-(\f n2-1)=\f n2+1>r$. Set
$G=rK_2\cup(k-r-1)K_1$. It is easily seen that $G$ is a graph of
order $k+r-1$ with $\a(G)=k-1$. For any $n$ vertices in $G$, the
corresponding induced subgraph by the $n$ vertices must be given by
$aK_2\cup bK_1$, where $a$ and $b$ are nonnegative integers such
that $2a+b=n$ and $a\le r$. Thus, $G$ is an $(n,r)$ graph of order
$k+r-1$. On the other hand, as $\a(G)=k-1$, $G$ has no independent
sets with $k$ vertices. Hence, $R(n,n(n-1)/2-r;k,1)>k+r-1$. Since
$k+r>n$, it follows from Theorem 3.1 that for any $(n,r)$ graph $G'$
of order $k+r$, $\a(G')\ge k+r-r=k$. Hence $R(n,n(n-1)/2-r;k,1)\le
k+r$ and so $R(n,n(n-1)/2-r;k,1)=k+r=\t{max}\{n,k+r\}$.
\par Now assume $r\geq \f {n-1}2.$ We first suppose $r\leq 2n-3k+2$.
Let $G$ be a graph of order $n$ with $e(G)\le r$. Taking $p=n$ and
$t=n-r$ in Theorem 3.2 we see that
$$\align\alpha(G)\geq \Big[ \f{n-[\f{n+4-2(n-r)}3]}2\Big]+1
=\Big[ \f{n-[\f{2r+4-n}3]}2\Big]+1\geq \Big[
\f{n-[\f{2(2n-3k+2)+4-n}3]}2\Big]+1=k.\endalign$$ Thus,
$R(n,\f{n(n-1)}2-r;k,1)\leq n.$ On the other hand, clearly
$R(n,\f{n(n-1)}2-r;k,1)>n-1.$ So $R(n,\f{n(n-1)}2-r;k,1)=n.$ To see
the result, we note that
$$2k-2+\big[\f{2r+4-n}3\big]\leq 2k-2+\big[\f{2(2n-3k+2)+4-n}3\big]=n.$$
Now we suppose $r>2n-3k+2.$ In this case, we have
$$2k-2+\Big[\f{2r+4-n}3\Big]\ge
 2k-2+\Big[\f{2(2n-3k+2)+4-n}3\Big]=n$$
and
$$\align 0\leq \f{2r+1-n}3&=\f{n+1-2(n-r)}3<\f{3k-2-(n-r)+1-2(n-r)}3\\&=k-\f
13-(n-r)\leq k-\f 13-2<k-1\endalign$$ and so $$ 0\leq
\big[\f{2r+1-n}3\big]\leq \f{2r+1-n}3<k-1.$$ Set
$$G=\big[\f{2r+1-n}3\big]K_3\bigcup
\big(k-1-\big[\f{2r+1-n}3\big]\big)K_2.$$ We claim that $G$ is an
$(n,r)$ graph. Clearly any induced subgraph of $G$ by $n$ vertices
can be written as $xK_1\cup yK_2\cup zK_3,$ where $x,y,z$ are
nonnegative integers satisfying $x+2y+3z=n$ and
$z\leq[\f{2r+1-n}3].$ If $x+y\leq n-r-1,$ then
$$\align\f{n+1-2(n-r)}3&\geq \big [\f{2r+1-n}3\big ]\geq
z=\f{n-2(x+y)+x}3\\&\geq\f{n-2(n-r-1)}3=\f{n+2-2(n-r)}3.\endalign$$
This is a contradiction. Hence $x+y\geq n-r$ and so
$$e(xK_1\cup yK_2\cup zK_3)=y+3z=n-(x+y)\leq r.$$
This shows that $G$ is an $(n,r)$ graph. It is evident that
$$\alpha(G)=\big[\f{2r+1-n}3\big]+k-1-\big[\f{2r+1-n}3\big]=k-1$$
 and
 $$|V(G)|=3\Big[\f{2r+1-n}3\Big]+2\Big(k-1-\Big[\f{2r+1-n}3\Big]\Big)
 =2k-2+\Big[\f{2r+1-n}3\Big].$$
 As $G$ has no independent sets with $k$ vertices and $G$ does
not contain any subgraphs with $n$ vertices and at least $r+1$
edges, we must have
$$R\Big(n,\f{n(n-1)}2-r;k,1\Big)>|V(G)|=2k-2+\big[\f{2r+1-n}3\big].$$
On the other hand, if $G$ is an $(n,r)$ graph of order
$2k-2+\big[\f{2r+4-n}3\big],$ by Theorem 3.2 we have
$$\alpha(G)\geq\Big[\f{2k-2+\big[\f{2r+4-n}3\big]-\big[\f{n+4-2(n-r)}3\big]}2\Big]
+1=k.$$ Hence $$R\Big(n,\f{n(n-1)}2-r;k,1\Big)\leq
2k-2+\big[\f{2r+4-n}3\big]$$ and so
$$R\Big(n,\f{n(n-1)}2-r;k,1\Big)=2k-2+\Big[\f{2r+4-n}3\Big]=\t{max}
\Big\{n, 2k-2+\Big[\f{2r+4-n}3\Big]\Big\}.$$
\par Summarizing the above we prove the theorem.
\pro{Corollary 4.1} Let $k$ and $n$ be positive integers with $k\ge
2$ and $n\ge 4$. Then
$$R(n,n(n-1)/2-n+2;k,1)
=\cases n&\t{if $n\ge 3k-4$,}
\\2k-2+[\f n3]&\t{if $n<3k-4$.}\endcases$$
\endpro
Proof. Putting $r=n-2$ in Theorem 4.1 we obtain the result.

\pro{Definition 4.1} Let $k,s\in\Bbb N$ with $s\le \b k2$. Let $G$
be a graph with at least $k$ vertices. If every subgraph of $G$
induced by $k$ vertices has at least $s$ edges, we say that $G$
satisfies the $(k,s)$ condition.\endpro
 \pro{Lemma 4.1} Let
$p,n,k,r,s\in\Bbb{N}$ with $n,k\geq 2,\ r<\b n2,\ s<\b k2$ and
$p\geq max\{n,k\}$. Let $G$ be a $(n,\b n2-r)$ graph of order $p$
satisfying the $(k,s)$ condition. Then for any $v\in V(G)$,
$$p-R(n,r;k-1,s)\leq d(v)\leq R(n-1,r;k,s)-1.$$\endpro
Proof. Let $x$ be a vertex in $G$ such that $d(x)=\Delta(G).$ If
$\Delta(G)\geq R(n-1,r;k,s),$ then $\Delta(G)\ge n-1$ and there are
$n-1$ vertices $x_1,\ldots,x_{n-1}\in\Gamma(x)$ such that
$e(G[\{x_1,$ $\ldots, x_{n-1}\}])>\binom {n-1}2-r.$ As
$\binom{n-1}2-r+n-1=\binom n2-r,$ we see that $e(G[\{x,x_1,\ldots,$
$x_{n-1}\}])>\binom n2-r.$ This contradicts the assumption. Hence
$\Delta(G)<R(n-1,r;k,s)$ and so $d(v)\leq \Delta(G)\leq
R(n-1,r;k,s)-1$ for any vertex $v$ in $G$.
\par Let $y$ be a vertex in $G$ such that $d(y)=\delta(G)$, and let
$V_y=V(G)-\{y\}\cup \Gamma(y)$. Then clearly
$|V_y|=p-1-d(y)=p-1-\delta(G).$ If $\delta(G)\ge p-R(n,r;k-1,s)$,
then $d(v)\ge \delta(G)\ge p-R(n,r;k-1,s)$ for any $v\in V(G)$. Now
assume $\delta(G)\le p-R(n,r;k-1,s)-1$. As $R(n,r;k-1,s)\ge k-1$ we
see that $\delta(G)\le p-k$ and so $|V_y|=p-1-\delta(G)\ge k-1$. For
any $k-1$ vertices $v_1,\ldots,v_{k-1}$ in $V_y,$ by the assumption
we have
$e(G[\{v_1,\ldots,v_{k-1}\}])=e(G[\{y,v_1,\ldots,v_{k-1}\}])\geq s$.
As $G$ is a $(n,\binom n2-r)$ graph, we see that $G[V_y]$ is also a
$(n,\binom n2-r)$ graph. Hence, by Definition 1.1 we have
$p-1-\delta(G)<R(n,r;k-1,s).$ Therefore, for any vertex $v$ in $G$,
$d(v)\geq \delta(G)\geq p-R(n,r;k-1,s),$ which completes the proof.
\pro{Theorem 4.2} Let $n,r,k,s$ be positive integers with $n,k\ge
3$. Then
$$R(n,r;k,s)\le R(n-1,r;k,s)+R(n,r;k-1,s).$$ Moreover, the strict
inequality holds when both $R(n-1,r;k,s)$ and $R(n,r;k-1,s)$ are
even.\endpro
 Proof. Clearly $R(n-1,r;k,s)+R(n,r;k-1,s)\ge \t{max}\
\{n-1+n,k+k-1\}>\t{max}\ \{n,k\}$. Thus, if $R(n,r;k,s)=\t{max}\
\{n,k\}$, the result is true. Now assume $R(n,r;k,s)>\t{max}\
\{n,k\}$. By the definition of $R(n,r;k,s),$ there is a $(n,\binom
n2-r)$ graph $G$ of order $R(n,r;k,s)-1$ satisfying the $(k,s)$
condition. For any vertex $v$ in $G$, by Lemma 4.1 we have
$$R(n,r;k,s)-1-R(n,r;k-1,s)\leq d(v)\leq R(n-1,r;k,s)-1.$$
It then follows that
$$R(n,r;k,s)\leq R(n-1,r;k,s)+R(n,r;k-1,s).$$
Moreover, the equality holds if and only if $G$ is a regular graph
with degree $R(n-1,r;k,s)-1$. If the equality holds, we must have
$$\align 2e(G)&=\big(R(n,r;k,s)-1\big)\big(R(n-1,r;k,s)-1\big)
\\&=\big(R(n-1,r;k,s)+R(n,r;k-1,s)-1\big)\big(R(n-1,r;k,s)-1\big).\endalign$$
This shows that either $R(n-1,r;k,s)$ or $R(n,r;k-1,s)$ is odd.
Hence, if  $R(n-1,r;k,s)$ and $R(n,r;k-1,s)$ are even, we have
$R(n,r;k,s)<R(n-1,r;k,s)+R(n,r;k-1,s)$, which completes the proof.
\par\q
\par Theorem 4.2 can be viewed as a generalization of the classical inequality ([1])
$R(n,k)\le R(n-1,k)+R(n,k-1)$ for $n,k\ge 3$.

\pro{Theorem 4.3} Let $k,n,r\in\Bbb N$ with $r\le n-1$. Then
\par $(\t{\rm i})$ $R(n,r;k,1)\ge R(n,r;k-1,1)+n-1.$
\par $(\t{\rm ii})$ $R(n,r;k,1)\ge (n-1)(k-1)+1.$
\endpro Proof. We first consider (i). By
the definition of $R(n,r;k-1,1)$, there is a $(k-1,\b{k-1}2-1)$
graph $G$ of order $R(n,r;k-1,1)-1$ satisfying the $(n,r)$
condition. Now we construct a new graph $G'$ by adding $n-1$ new
vertices to $G$ and joining every new vertex and each vertex in $G$.
Since $G$ does not contain a copy of $K_{k-1}$, we see that $G'$
does not contain a copy of $K_k$. For fixed
$s\in\{1,2,\ldots,n-1\}$, if we choose $n$ vertices in $G'$
containing exactly $s$ vertices in $G$, then the subgraph of $G'$
induced by the $n$ vertices has at least $s(n-s)$ edges. As
$s(n-s)\ge n-1\ge r$ and $G$ satisfies the $(n,r)$ condition, we see
that $G'$ also satisfies the $(n,r)$ condition. Note that the order
of $G'$ is $R(n,r;k-1,1)-1+n-1$. We then have
$R(n,r;k,1)>R(n,r;k-1,1)-1+n-1$. This proves (i).
\par Now let us consider (ii). It is clear that $R(n,r;2,1)=n$.
Thus, using (i) we see that
$$\align R(n,r;k,1)&=R(n,r;2,1)+\sum_{s=3}^k\big(R(n,r;s,1)-R(n,r;s-1,1)\big)
\\&\ge n+(k-2)(n-1)=(n-1)(k-1)+1.\endalign$$
 This proves (ii) and hence the
theorem is proved.
\section*{5. The upper
bound for $R(4,3;k,1)$}
\par Let $k\in\{2,3,\ldots\}$. By Definition 1.1, $R(4,3;k,1)$ is the smallest positive integer $p$ such
 that for any graph $G$ of order $p$, either $G$
has a subgraph induced by $4$ vertices with at least $4$ edges, or
$G$ contains an independent set with $k$ vertices. Every subgraph of
$(k-1)K_3$ induced by $4$ vertices has at most three edges and the
independence number of $(k-1)K_3$ is $k-1$. Thus
$R(4,3;k,1)>3(k-1)$.
 \pro{Theorem 5.1} For $k=3,4,5,6$ we have $R(4,3;k,1)=3k-2$.\endpro
Proof. By the previous argument or Theorem 4.3(ii) we have
$R(4,3;k,1)\ge 3k-2$. Clearly $R(4,3;2,1)=4$ and $R(3,3;3,1)=3$.
Thus, using Theorem 4.2 we see that
$$R(4,3;3,1)\le R(4,3;2,1)+R(3,3;3,1)=4+3=7.$$
Hence $R(4,3;3,1)=7$.
\par For $k\in\{4,5,6\}$ let $G$ be a $(4,3)$ graph of order $3k-2$.
 If $\a(G)<k$, for any vertex $v$ in
$G$, by Lemma 4.1 we have
$$3k-2-R(4,3;k-1,1)\le d(v)\le R(3,3;k,1)-1=k-1.\tag 5.1$$
For a $(4,3)$ graph $G$ of order $10$ with $\a(G)<4$, it follows
from (5.1) that $G$ is a 3-regular graph with $g(G)\ge 5$. As
$R(3,4)=9$, we must have $\a(G)\ge 4$. This is a contradiction.
Hence for any $(4,3)$ graph $G$ of order $10$ we have $\a(G)\ge 4$
and so $R(4,3;4,1)=10$.
\par Now let $G$ be a $(4,3)$ graph of order $13
$. If $\a(G)<5$, from (5.1) we see that $3\le d(v)\le 4$ for any
$v\in V(G)$. If $d(v)=\Delta(G)=4$, $\Gamma(v)=\{v_1,v_2,v_3,v_4\}$
and $x_1,x_2\in\Gamma(v_1)-\{v\}$, then clearly
$\{x_1,x_2,v_2,v_3,v_4\}$ is an independent set of $G$ and so
$\a(G)\ge 5$. If $G$ is a 3-regular $(4,3)$ graph of order $13$,
then $g(G)\ge 5$ and for given $v\in V(G)$ there are three vertices
$u_i(i=1,2,3)$ in $G$ such that $d(u_i,v)\ge 3\ (i=1,2,3)$. As
$g(G)\ge 5$, we may assume that $u_i$ is not adjacent to $u_j$, then
clearly $\{u_i,u_j\}\cup\Gamma(v)$ is an independent set of $G$ and
so $\a(G)\ge 5$. By the above, for any $(4,3)$ graph $G$ of order
$13$ we have $\a(G)\ge 5$ and so $R(4,3;5,1)=13$.
\par Suppose
that $G$ is a $(4,3)$ graph of order $16$. If $\a(G)<6$, from (5.1)
we see that  $3\le d(v)\le 5$ for any $v\in V(G)$. If
$d(v)=\Delta(G)=5$, $\Gamma(v)=\{v_1,\ldots,v_5\}$ and
$x_1,x_2\in\Gamma(v_1)-\{v\}$, then clearly
$\{x_1,x_2,v_2,v_3,v_4,v_5\}$ is an independent set of $G$ and so
$\a(G)\ge 6$. If $d(v)=\Delta(G)=4$, $\Gamma(v)=\{v_1,\ldots,v_4\}$
and $\Gamma(v_1)=\{v,x_1,x_2,x_3\}$, then clearly
$\{x_1,x_2,x_3,v_2,v_3,v_4\}$ is an independent set of $G$ and so
$\a(G)\ge 6$. If $d(v)=\Delta(G)=4$, $\Gamma(v)=\{v_1,\ldots,v_4\}$
and $d(v_1)=\cdots=d(v_4)=3$, there are three vertices
$u_i(i=1,2,3)$ in $G$ such that $d(u_i,v)\ge 3\ (i=1,2,3)$. As
$g(G)\ge 5$, we may assume that $u_i$ is not adjacent to $u_j$. Then
clearly $\{u_i,u_j,v_1,\ldots,v_4\}$ is an independent set of $G$
and so $\a(G)\ge 6$. If $d(v)=\Delta(G)=3$ and
$\Gamma(v)=\{v_1,v_2,v_3\}$, then $G$ is 3-regular and there are six
vertices $u_i(i=1,2,\ldots,6)$ in $G$ such that $d(u_i,v)\ge 3\
(i=1,2,\ldots,6)$. As $g(G)\ge 5$ and $R(3,3)=6$, there are three
vertices $u_i,u_j,u_k\in\{u_1,\ldots,u_6\}$ such that
$\{u_i,u_j,u_k\}$ is an independent set and so
$\{u_i,u_j,u_k,v_1,v_2,v_3\}$ is also an independent set of $G$.
This shows that $\a(G)\ge 6$. Hence, by the above, for any $(4,3)$
graph $G$ of order $16$ we have $\a(G)\ge 6$ and so $R(4,3;6,1)=16$.
This completes the proof.

 \pro{Theorem 5.2} For $k=3,4,5,\ldots$ we have
$$R(4,3;k,1)\le R(4,3;k-1,1)+\f 32+\sqrt{R(4,3;k-1,1)-\f 34}.$$
\endpro
Proof. By the definition of $R(4,3;k,1),$ there is a $(4,3)$ graph
$G$ of order $R(4,3;k,1)-1$ with $\a(G)<k$. For any vertex $v$ in
$G$, by Lemma 4.1 we have
$$R(4,3;k,1)-1-R(4,3;k-1,1)\leq d(v)\leq R(3,3;k,1)-1=k-1.\tag 5.2$$
If $\delta(G)\le 2$, as $R(4,3;k,1)\ge R(4,3;3,1)=7$ we have
$\delta(G)\le \sqrt {R(4,3;k,1)-2}$. If $\delta(G)\ge 3$, then
clearly $g(G)\ge 5$ and so $\delta(G)\le \sqrt {R(4,3;k,1)-2}$ by
[1, p.105] or [10, Proposition 1]. Now, from the above we deduce
that
$$R(4,3;k,1)-1-R(4,3;k-1,1)\le \delta(G)\le \sqrt{R(4,3;k,1)-2}.
\tag 5.3$$ Set $x_n=R(4,3;n,1)$. Then $(x_k-x_{k-1}-1)^2\le x_k-2$
and so $x_k^2-(2x_{k-1}+3)x_k+(x_{k-1}^2+2x_{k-1}+3)\le 0.$ This
yields
$$\align x_k\le&
\f{2x_{k-1}+3+\sqrt{(2x_{k-1}+3)^2-4(x_{k-1}^2+2x_{k-1}+3)}}2
\\&=x_{k-1}+\f 32+\sqrt{x_{k-1}-\f 34}.\endalign$$
 This completes the proof.
\pro{Theorem 5.3} Let $0<\ep\le 1$ and $k\in\Bbb N$ with $k\ge 6$.
Then
$$R(4,3;k,1)\le \f{(k-6)(k+6+2a)}{4-\ep}+16<\f{(k+a)^2}{4-\ep}$$
and $$R(4,3;k,1)-R(4,3;k-1,1)<1+\f{k+a}{\sqrt{4-\ep}},$$
 where
$$a=\f{5-1.5\ep}{2-\sqrt{4-\ep}}-6.$$\endpro
Proof. As
$$a+6=\f{5-1.5\ep}{2-\sqrt{4-\ep}}\ge\f{5-1.5}{2-\sqrt
3}>8>4\sqrt{4-\ep},$$ we have $\f{(a+6)^2}{4-\ep}>16$. Now we prove
the first part by induction on $k$. Since $R(4,3;6,1)=16$, the
result is true for $k=6$. Suppose $k\ge 7$ and
$$R(4,3;k-1,1)\le
\f{(k-1-6)(k-1+6+2a)}{4-\ep}+16=\f{(k-1+a)^2}{4-\ep}
-\f{(6+a)^2}{4-\ep}+16.$$ Then
$$\align &R(4,3;k-1,1)+\f 32+\sqrt{R(4,3;k-1,1)-\f 34}
\\&<\f{(k-1+a)^2}{4-\ep}-\f{(6+a)^2}{4-\ep}+16+\f 32+\f{k-1+a}{\sqrt{4-\ep}}
\\&=\f{(k+a)^2}{4-\ep}-\f{(6+a)^2}{4-\ep}+16+\f 32-\Big(\f
2{4-\ep}-\f 1{\sqrt{4-\ep}}\Big)(k+a-1)-\f 1{4-\ep}.\endalign$$ As
$k+a-1\ge a+6=\f{5-1.5\ep}{2-\sqrt{4-\ep}}$, we see that
 $(k+a-1)(2-\sqrt{4-\ep})\ge 5-1.5\ep=\f 32(4-\ep)-1$ and so
$$\Big(\f 2{4-\ep}-\f 1{\sqrt{4-\ep}}\Big)(k+a-1)+\f 1{4-\ep}>\f
32.$$ Hence, by the above and Theorem 5.2 we obtain
$$\align R(4,3;k,1)&\le  R(4,3;k-1,1)+\f 32+\sqrt{R(4,3;k-1,1)-\f 34}
\\&<\f{(k+a)^2}{4-\ep}-\f{(6+a)^2}{4-\ep}+16= \f{(k-6)(k+6+2a)}{4-\ep}+16<\f{(k+a)^2}{4-\ep}.\endalign$$
 From the above and (5.3) we also deduce that
$$\align &R(4,3;k,1)-R(4,3;k-1,1)-1\\&\le\sqrt{R(4,3;k,1)-2}
<\sqrt{\f{(k+a)^2}{4-\ep}-\f{(6+a)^2}{4-\ep}+14}<\f{k+a}{\sqrt{4-\ep}}.\endalign$$
This completes the proof.
 \pro{Theorem 5.4} Let $0<\ep\le 1$ and
$a=\f{5-1.5\ep}{2-\sqrt{4-\ep}}-6.$ Let $G$ be a graph of order
$p\ge 16$ with $g(G)\ge 5$. Then
$$\a(G)\ge [\sqrt{(4-\ep)(p-16)+(a+6)^2}-a].$$
\endpro
Proof. Set $k=[\sqrt{(4-\ep)(p-16)+(a+6)^2}-a].$ Then clearly $k\ge
6$. Using Theorem 5.3 we see that
$$p\ge \f{(k+a)^2}{4-\ep}-\f{(6+a)^2}{4-\ep}+16\ge R(4,3;k,1).$$
Since $g(G)\ge 5$, $G$ must be a $(4,3)$ graph. Hence $\a(G)\ge k$.
This proves the theorem.
\section*{6. Some open conjectures on
Ramsey numbers}
\par In this section we risk to pose the following conjectures on
Ramsey numbers.
 \newline{\bf Conjecture 6.1.} For any positive integer $n\ge 2$ we have
$$\f{n-1}{R(3,n)-1}>\f n{R(3,n+1)-1}\qtq{and
so}R(3,n+1)>\f{nR(3,n)-1}{n-1}.$$
\par As $\f 12>\f 25>\f 38>\f 4{13}>\f 5{17}>\f 6{22}>\f 7{27}>\f
8{35}$, from (1.1) we know that Conjecture 6.1 is true for
$n\in\{2,3,\ldots,8\}$. If the conjecture is true, we have
$R(3,10)>\f{9R(3,9)-1}8>40$. It is now known ([8]) that $40\le
R(3,10)\le 43$.
\newline{\bf Conjecture 6.2.} Let $\{L_n\}$ be the Lucas sequence
defined by $L_0=2,\ L_1=1$ and $L_{n+1}=L_n+L_{n-1}(n\ge 1)$. For
$k=3,4,5,\ldots$ we have $R(k,k)=4L_{2k-5}+2.$
\par Conjecture 6.2 is true for $k=3,4$. By
Conjecture 6.2, we have $R(5,5)=46$, $R(6,6)=118$ and $R(7,7)=306$.
 Since $L_{2(n+1)}=3L_{2n}-L_{2(n-1)}$, Conjecture 6.2 is
equivalent to
$$R(k,k)=3R(k-1,k-1)-R(k-2,k-2)-2\qtq{for}k\ge 3.\tag 6.1$$
It is well known that
$$L_n=\Big(\f {1+\sqrt 5}2\Big)^n+ \Big(\f
{1-\sqrt 5}2\Big)^n.$$ Thus, by Conjecture 6.2,
$$\aligned R(k,k)&=4\Big\{\Big(\f {1+\sqrt 5}2\Big)^{2k-5}- \Big(\f
{\sqrt 5-1}2\Big)^{2k-5}\Big\}+2 \\&= 128\Big\{\Big(\f {3+\sqrt
5}2\Big)^{k}- \Big(\f {3-\sqrt 5}2\Big)^{k}\Big\}+2 .\endaligned$$
Hence,
$$R(k,k)\sim 128\Big(\f{3+\sqrt 5}2\Big)^k\qtq{as}k\rightarrow
+\infty.\tag 6.2$$ We note that $\f{3+\sqrt 5}2\approx 2.618$. It is
known that $(\sqrt 2)^k< R(k,k)\le 4^k.$ P. Erd\H{o}s offered
$\$350$ to ask the value of $\lim\limits_{k\rightarrow
\infty}R(k,k)^{\f 1k}$ (see [3, p.10]). If the limit exists, it
should be $\f{3+\sqrt 5}2$ by Conjecture 6.2.
\newline{\bf Conjecture 6.3.}
For $n=2,3,4,\ldots$ we have
$$\sum_{r=1}^{n(n-1)/2}R(n,r;3,1)=R\big(3,\f {n(n+1)}2-1\big).$$
\par Conjecture 6.3 is true for $n=2,3,4$. Since
$$\sum_{r=1}^{10}R(5,r;3,1)=14+11+9+9+7+7+5+5+5+5=77,$$
by Conjecture 6.3 we have $R(3,14)=77$. It is known ([8]) that
$66\le R(3,14)\le 78$.
\newline{\bf Conjecture 6.4.} For $k=1,2,3,\ldots$
we have $R(4,3;k,1)=3k-2.$
\par From Theorem 5.1 we know that Conjecture 6.4 is true for $k\le 6$.

\end{document}